\documentclass{article}

\usepackage{amsmath}
\usepackage{amsthm}

\usepackage{subfigure}
\usepackage{graphicx}
\usepackage{epstopdf}
\usepackage{fancybox}

\newcommand{\Vecb}[2]{\mathbf{#1}_{#2}}

\newcommand{\Matb}[2]{\mathbf{#1}_{#2}}

\newcommand{\state}[1]{\Vecb{x}{#1}}
\newcommand{\obs}[1]{\Vecb{z}{#1}}

\newcommand{\xEst}[2]{\hat{\mathbf{x}}_{#1|#2}}
\newcommand{\xErr}[2]{\tilde{\mathbf{x}}_{#1|#2}}
\newcommand{\xErrT}[2]{\tilde{\mathbf{x}}_{#1|#2}^\top}
\newcommand{\covEst}[2]{\mathbf{P}_{#1|#2}}

\newcommand{\pNoiseb}[1]{\Vecb{v}{#1}}
\newcommand{\oNoiseb}[1]{\Vecb{w}{#1}}

\newcommand{\tc}{k}
\newcommand{\tl}{k-1}

\newcommand{\E}[1]{\mathrm{E}\left[#1\right]}
\begin{document}


\title{{\bf Gaussianity and the Kalman Filter:\\
       A Simple Yet Complicated Relationship}\footnote{Published in {\em Journal de Ciencia e Ingeniería}, vol.\ 14, no.\ 1, pp.\ 21-26, 2022.}}
\author{{\bf Jeffrey K.\ Uhlmann}\\
Dept. of Electrical Engineering and Computer Science\\
University of Missouri-Columbia\\~\\
{\bf Simon J.\ Julier}\\
Department of Computer Science\\ University College London}
\date{}     
   
\maketitle
\thispagestyle{empty}

\begin{abstract}
  One of the most common misconceptions made about the Kalman filter when
  applied to linear systems is that it requires an assumption that all error and
  noise processes are Gaussian. This misconception has frequently led to
  the Kalman filter being dismissed in favor of complicated and/or purely
  heuristic approaches that are supposedly ``more general'' in that they
  can be applied to problems involving non-Gaussian noise. The fact is that
  the Kalman filter provides rigorous and optimal performance guarantees
  that do not rely on any distribution assumptions beyond mean and error
  covariance information. These guarantees even apply to use of the Kalman 
  update formula when applied with nonlinear models, as long as its other
  required assumptions are satisfied. Here we discuss misconceptions about 
  its generality that are often found and reinforced in the literature, especially
  outside the traditional fields of estimation and control.

\end{abstract}

\noindent {\em Thematic areas}: Educational Engineering; Engineering History.

\maketitle

\section{Introduction}

The Kalman filter~\cite{article:kalman60} is among the most versatile and
widely-used tools in engineering. More than 50 years after Kalman published
his original paper, his work still inspires hundreds of papers each year. Some
of these papers explore new applications of the algorithm in approaches which
range from industrial process control, to robotics, and even to meta-analysis of
election forecasts. 

Unfortunately, the full potential of the Kalman filter is often not appreciated and exploited due to
misconceptions which have persisted since the earliest days after its emergence
as a critical component in aeronautical and space applications in the 1960s.
Among the most common misconceptions is that the Kalman filter can only be
rigorously derived and applied to linear systems in which all the error and
noise processes are \emph{strictly\/} Gaussian. More specifically, it is
commonly believed --- and frequently stated implicitly or explicitly --- that
the use of a Kalman filter in the presence of non-Gaussian error processes is at the
very least a sub-optimal heuristic approach that may perform well in practice if
errors are approximately Gaussian but that it is mathematically non-rigorous and
cannot be expected to perform well if the errors are strongly non-Gaussian.
Figure~\ref{fig:deluded_quotes} gives evidence of this in form of a 
quotes sampled from the literature across a range of fields
including machine learning, estimation, systems engineering and end-user
applications. (These quotes should not be taken as undermining the integrity
of the sources from which they are taken but rather as examples of how
someone new to the study of linear and estimation might come to believe
that the Kalman filter requires Gaussianity assumptions.)

\begin{figure}

\begin{center}
     \fbox{%
       \begin{minipage}{0.95\linewidth}\small
\begin{itemize}

\item
``{\em Since the Kalman framework requires Gaussian
distributions, the model can only be constructed if
...~}''~\cite{inproceedings:yardim07}


\item
``{\em The Kalman filter which is used for integrated navigation
requires Gaussian variables ... a multimodal un-symmetric 
distribution has to be approximated with a
Gaussian distribution before being used in the Kalman
filter.}''~\cite{article:schon05}


\item
``{\em ...can be best reconciled with the KF (which requires Gaussian
probability distributions) by making the assumption that
...~}''~\cite{article:osborn10}


\item
`` {\em [The] Kalman filter requires Gaussian prior $f(x_0)$
...~}''~\cite{article:suzdaleva08}


\item
``{\em Notice that each of the distributions can be effectively approximated by a Gaussian. This is a very important result for the
operation for many systems, especially the ones based on a Kalman filter since the filter explicitly requires Gaussian distributed
noise on measurements for proper operation.}''~\cite{techreport:aksoy14}


\item
{\em ...importance sampling ... relaxes the assumption of Gaussian observation errors required by the basic Kalman filter.}
~\cite{article:knape11}

\item
``{\em However, the KF requires Gaussian initial conditions, therefore
...~}''\cite{inproceedings:shaferman09}


\item
``{\em Kalman filters are Bayes filters that represent posteriors with Gaussians... Kalman filter mapping relies on three basic assumptions ... Gaussian noise... the initial uncertainty must be Gaussian.}''~\cite{Thrun:2003:RMS:779343.779345}


\item
``{\em The Kalman filter is a very efficient optimal filter; however
it has the precondition that the noises of the process
and of the measurement are Gaussian ... when the measurement is not a Gaussian
distribution, the Kalman filter cannot be used.}''~\cite{inproceedings:rosenberg97}


\item
``{\em The Kalman filter assumes that the posterior density at every time step is
Gaussian and hence exactly and completely parameterized by two parameters, its
mean and covariance.}''


\item
``{\em Kalman filter cannot be used here for inference because the measurement
does not involve additive Gaussian noise.}''~\cite{inproceedings:pfister09}


\item
``{\em The KF requires models defined by linear Gaussian probability density functions.}''~\cite{0967-3334-34-7-781}


\item ``{\em The most attractive advantage of the Kalman filter
    lies in its optimal estimation in the sense of minimum mean
    squared prediction errors. However, the optimality of the Kalman
    filter requires two restrictive prerequisites, linear state-space
    models and independent Gaussian white noise for both process and
    measurements.}''~\cite{Bai200565}

  
\end{itemize}       \end{minipage}}
\end{center}

\caption{Quotes in the literature that could be interpreted as suggesting that the Kalman filter can only be
applied when errors are Gaussian-distributed.}
\label{fig:deluded_quotes}

\end{figure}

An unfortunate consequence of such misconceptions is that it is common for the
Kalman filter to be dismissed from consideration for applications simply because
the errors are known to be non-Gaussian. In this paper, we discuss how the 
Kalman filter can be derived as the optimal solution to the filtering problem 
under differing sets of assumptions, and we emphasize that the assumptions
required for a particular derivation are not necessarily required in general
for the optimality of the filter. In particular, we emphasize that linearity
does not imply Gaussianity; minimizing mean-squared error does not imply
Gaussianity; and that the Kalman Filter is MMSE-optimal without any assumptions
of Gaussianity.

In the next section we briefly summarize the linear
estimation/filtering problem and historically how the same optimal
solution has been derived from two very different perspectives with
very different assumptions.In particular, we emphasize that the Kalman
filter can be applied rigorously -- and optimally -- to systems with
errors from any probability distribution with finite first and second
moments

\section{The Estimation Problem: Kalman vs. Bayes}
\label{problemStatement}

\subsection{Estimation Problem}

Consider a linear system of the form
\begin{equation}
\state{\tc}=\Matb{F}{\tl}\state{\tl}+\pNoiseb{\tl},
\label{eqn:prediction_model}
\end{equation}
where $\state{\tc}$ is the state at timestep $\tc$, $\Matb{F}{\tl}$ is
the state transition matrix, and $\pNoiseb{\tl}$ is the additive
process noise. We assume that this noise is independent from timestep, 
is zero mean, and its covariance is known. However, we do
\emph{not\/} assume that it is Gaussian-distributed.

The observation model for a sensor measurement of the state of the system has
the same linear form:
\begin{equation}
\obs{\tc}=\Matb{H}{\tc}\state{\tc}+\oNoiseb{\tc},
\end{equation}
where where $\Matb{H}{\tc}$ is a linear transformation from 
system to observation coordinates and $\oNoiseb{\tc}$ is the
observation noise. Similar to $\pNoiseb{\tl}$, it is assumed to be
zero-mean, independent, but does not have to be Gaussian.

The filter is initialized at time step 0 with an estimate $\xEst{0}{0}$. The
error in this estimate is zero mean and has a covariance $\covEst{0}{0}$. Again,
this does not have to be Gaussian.

Given an initial condition and a sequence of measurements, the goal is to
compute an estimate $\xEst{i}{j}$ at timestep $i$ based on all observations up
to timestep $j$, along with its associated error covariance $\covEst{i}{j}$. The
relationship between the estimate and the state is given by
\begin{equation}
\xEst{i}{j}=\state{i}+\xErr{i}{j},
\end{equation}
where $\xErr{i}{j}$ is the error. The mean squared error in this
estimate is
\begin{equation}
   \covEst{i}{j} ~=~ \E{\xErr{i}{j}\xErrT{i}{j}}.
\end{equation} 
Given that the covariance of any error distribution can never 
be determined {\em exactly} in practice, the above equality
can be relaxed to the more conservative and 
practically-achievable requirement:
\begin{equation}
   \covEst{i}{j} ~\geq~ \E{\xErr{i}{j}\xErrT{i}{j}},
\end{equation}
where the covariance matrix $\covEst{i}{j}$ can now be
interpreted as representing the best available upper bound
on the expected squared error associated with $\xEst{i}{j}$. 
As we discuss later, the practical question is whether an estimate
which obeys this guarantee is sufficient for the problem at hand.

In Kalman's original paper he derived his now-eponymous filter from
the perspective of $\ell_2$-norm error minimization via othogonal
projections. He also noted (his Corollary~1) that if all errors are
assumed Gaussian then the system mean-and-covariance estimate
can be interpreted as parameterizing a Gaussian distribution that 
represents the exact error distribution conditioned on the sequence 
of observations. In more contemporary parlance, Kalman derived
a minimum-mean-squared error (MMSE) optimal filter. What is critical 
to note, however, is that the MMSE optimality of the filter does not 
in any way depend on such an assumption.  

While MMSE optimality can be obtained with broad generality, i.e., 
with relatively weak assumptions, it does so at the expense of a
strong probabilistic interpretation. For example, the Kalman filter
guarantees that the expected squared error of the system estimate
decreases at a certain rate, but it does not provide information 
necessary to answer a question about the probability that, e.g.,
its mean position estimate is within one meter of the true position.
This motivated Ho and Lee~\cite{Ho1964} to re-derive the optimal
filter from a Bayesian perspective by replacing the mean-covariance
pair with the full probability density of the state, 
$f_{k|k}(\state{k}|\obs{1:k})$. The prediction is then determined
by the Chapman-Kolmogorov Equation,
\begin{equation}
f_{\tc|\tl}(\state{\tc}|\obs{1:\tl})=\int{f(\state{\tc}|\mathbf{x}^\prime)
f_{\tl|\tl}(\mathbf{x}^\prime|\obs{1:k})}\mathrm{d}\mathbf{x}^\prime
\end{equation}
where $f(\state{\tc}|\mathbf{x}^\prime)$ is the state transition density which
encodes the process model. The update is then given from Bayes rule,
\begin{equation}
f_{\tc|\tc}(\state{\tc}|\obs{1:\tc})=\frac{f(\obs{\tc}|\state{\tc})f_{\tc|\tl}(\state{\tc}|\obs{1:\tl})}
{f(\obs{1:\tc})},
\end{equation}
where $f(\obs{\tc}|\state{\tc})$ is the measurement likelihood model. This
incorporates the effects of the observation model.

From this Ho and Lee provided an alternate proof of Kalman's corollary
that the Kalman filter is Bayes-optimal with all mean and covariance 
estimates interpreted as parameters for Gaussian densities corresponding 
to assumed-Gaussian error processes. Under these assumptions a mean
and covariance estimate from the Kalman filter does not just represent 
the first two moments of an otherwise unknown probability distribution,
it can be interpreted as the {\em exact} uncertainty distribution for the
state. Having access to the full error distribution is clearly more informative,
but it comes at the cost of assumptions that cannot be realistically satisfied 
in almost any nontrivial real-world application. Even beyond Gaussianity,
the Bayesian derivation strictly requires {\em exact} and {\em
complete} knowledge about the full statistics of all noises. The least-squares
derivation, by contrast, easily accommodates conservative covariance estimates
and is thus much more consistent with the practical reality that the error
covariance of a sensor can never be ascertained with infinite precision.

While it is debatable whether or not the Bayesian interpretation is
pedagogically more accessible or intuitive than squared-error minimization,
there is no question that it appears more prominently in textbooks and
introductory expositions of the Kalman filter.  It should not be surprising,
therefore, that lingering concerns might exist (e.g., as suggested by the quotes
in Figure~\ref{fig:deluded_quotes}) that Gaussian-distributed noises are somehow
more conducive to good filter performance based on their special role in
Bayesian derivations.
It would seem that the MMSE derivation should be sufficient to dispel such
concerns, but its abstract mathematical guarantees apparantly cannot break the
intuitive link between covariances and Gaussians that is so firmly ingrained by
the Bayesian interpretation.

\section{Conclusions}
\label{sct:discussion}

The Procrustean application of Bayes' rule to
derive the Kalman filter may be suitable as a pedagogical exercise,
but care must be taken to ensure that the assumptions required for
the method of derivation are not confused with assumptions
that are required in general for effective use of the filter.
Linearity does not imply Gaussianity. Minimizing mean squared error
does not imply Gaussianity. The Kalman Filter is MMSE-optimal 
without any assumptions of Gaussianity. 
In this note we have attempted to highlight this fact so 
that the generality and optimality of the Kalman filter
can be more fully appreciated.

\end{document}